\documentclass[11pt]{article}
\newtheorem{theorem}{Theorem}
\newtheorem{lemma}{Lemma}
\newtheorem{corollary}{Corollary}

\usepackage{epsfig}

\title{Bells, Motels and Permutation Groups}
\author{Gary McGuire\\
School  of Mathematical Sciences\\
University College Dublin\\
Ireland}
\begin{document}
\maketitle

\begin{abstract}

This article is about the mathematics of ringing the changes. We
describe the mathematics which arises from a real-world activity,
that of ringing the changes on bells. We present Rankin's solution
of one of the famous old problems in the subject.
This article was written in 2003.
\end{abstract}

\section{Introduction: Motels and Bells}

In chapter 6 of his book ``Time Travel and other Mathematical
Bewilderments'' \cite{Gard}, Martin Gardner discusses the
following problem, a special case of D. H. Lehmer's  ``motel
problem'' from \cite{Lehmer}.

\bigskip

{\sl Mr.\ Smith manages a motel.
It consists of $n$ rooms in a straight row.
There is no vacancy.
Smith is a psychologist who plans to study the effects
of rearranging his guests in all possible ways.
Every morning he gives them a permutation.
The weather is miserable, raining almost daily.
To minimize his guests' discomfort, each daily rearrangement
is made by exchanging only the occupants of two adjoining rooms.
Is there a simple algorithm that will run through all
possible rearrangements by switching
only one pair of adjacent occupants
at each step}?

\bigskip

For the purposes of this article, we refer to this problem as the
{\it motel problem}. Before discussing a solution we consider
another problem, which at first might seem unrelated.  This
problem concerns the change-ringing of bells, so we first provide
a brief introduction to this topic and an explanation of how
permutations arise in the ringing of bells.

Around the year 1600 in England it was discovered that by altering
the fittings around each bell in a bell tower, it was possible for
each ringer to maintain precise control of when his (there were no
female ringers then) bell sounded.   This enabled the ringers to
ring the bells in any particular order, and either maintain that
order or change the order in a precise way.

Suppose there are $n$ bells being rung, numbered $1, 2, 3, \ldots
,n$ in order of pitch, number 1 being the highest.  When the bells
are rung in order of descending pitch $1, 2, 3, \ldots ,n$, we say
they are being rung in {\it rounds}. A change in the order of the
bells, such as rounds $1 \ 2 \ 3 \ 4 \ 5$ being changed to $2\ 1\
4\ 3\ 5$, can be considered as a permutation in the symmetric
group on five objects. In the years 1600--1650 a new craze emerged
where the ringers would continuously change the order of the bells
for as long as possible, while not repeating any particular order,
and return to rounds at the end. This game evolved into a
challenge to ring the bells in every possible order, without any
repeats, and return to rounds. We will give a precise statement of
this challenge shortly. However, the reader can now see why
bell-ringers were working with permutations. It has been pointed
out before \cite{WhiteSted} that permutations were first studied
in the 1600's in the context of the ringing of bells in a certain
order, and not in the 1770's by Lagrange, in the context of roots
of polynomials.  And, by the way, the craze continues to this day.

\bigskip

A solution to the motel problem is a sequence of all the
arrangements (orderings) of $n$ objects, with the property that
each arrangement is obtained from its predecessor by a single
interchange of two adjacent objects. Algorithms for generating the
$n!$ arrangements subject to this condition can be found in S. M.
Johnson \cite{Johnson} and H. F. Trotter \cite{Trotter}; see also
the papers by D. H. Lehmer \cite{Lehmer} and M. Hall Jr.--D. Knuth
\cite{HallK}. Johnson and Trotter discovered the same algorithm,
and we describe their algorithm below. Gardner \cite{Gard} credits
Steinhaus \cite{Steinhaus} as being the first to discover this
method. Such algorithms are similar to what bell-ringers are
trying to find, as we shall now explain.
In fact, the Johnson-Steinhaus-Trotter algorithm was discovered
in the 17th century by bell-ringers.
We elaborate on this in section \ref{jst}.

Let us refer to the $n$ bells listed in a particular order or
arrangement as a {\it row}.  For example, when $n=5$, rounds is
the row $1 \ 2 \ 3 \ 4 \ 5$. We will use cycle notation for
permutations. Permutations act on the positions of the bells, and
not on the bells themselves; for example, the permutation $(1\
2)(3\ 4)$ changes the row $5\ 4\ 3\ 2\ 1$ into the row $4\ 5\ 2\
3\ 1$. It is important not to confuse the numbers standing for the
bells with the numbers in the cycles, which refer to the positions
of the bells.  Permutations act on the right, i.e., $XY$ means
first do the permutation $X$, then do $Y$. As usual, $S_n$ and
$A_n$ denote the symmetric group and alternating group on $n$
objects.

The task enjoyed by bell-ringers, which we refer to as the {\it
ringers problem}, is to make a list of  rows subject to the
following rules:

\bigskip

(1) The first and last rows must be rounds.

\bigskip

(2) No row may be repeated (apart from rounds which appears twice,
first and last).

\bigskip

(3) Each bell may only change place by one position when moving
from one row to the subsequent row.

\bigskip

(4R)  No bell occupies the same position for more than two
successive rows.

\bigskip

A {\it method} is any set of rows that obeys these rules.  The
origins of rules (1) and (2) have been explained in the
introduction. Rule (3) is there for physical reasons. The purpose
of rule (4R) is mainly to make ringing interesting for the
ringers. For more details see \cite{Fletch} or \cite{WhiteRing}.
Clearly, $n!+1$ is the longest possible length that a list of rows
satisfying rules (1), (2), (3), and (4R) could have, and a list of
this length is called an {\it extent}.  The ringers problem is
essentially to construct methods of various lengths, and of
particular interest is the construction of extents.

The motel problem  is not identical to the ringers problem, for if
we replace rule (4R) by the following:

\bigskip

(4M)  Only two bells may change place when moving from one row to
the subsequent row

\bigskip then a list of $n!+1$ rows satisfying rules (1), (2), (3),
and (4M) is precisely a solution to the motel problem.

We can now see that there is certainly a  similarity between the
motel problem and the ringers problem. Some previous {\it Monthly}
articles have considered the motel problem \cite{HallK},
\cite{Lehmer}, and other {\it Monthly} articles have considered
the ringers problem, \cite{Dick}, \cite{Fletch}, \cite{WhiteRing},
\cite{WhiteSted}.
%but no article has pointed out the
%similarities between them.

We should mention that in practice ringers do have some other
rules, but rules (1), (2), (3), (4R) are the most important. We
ignore the other rules, and some other ringing matters, for the
purposes of this article. For more details see \cite{Fletch} or
\cite{WhiteRing}.

In the following sections we discuss these two problems and some
solutions.  In section \ref{jst} we present the first solution to
the motel problem, which was discovered by ringers in the 1600's.
Section \ref{plainhunt} discusses the simplest solution to the
ringers problem, and section \ref{plainbob} gives a more
complicated solution.  We shall see how ringers were manipulating
cosets of a dihedral subgroup of the symmetric group in an effort
to meet their challenge.  Section \ref{ccc} proposes a new method.
In section \ref{groupform} we present a general group-theoretic
framework for the problems. We discuss approaches to finding
extents in section \ref{leads}, and give some proofs of
impossibility.   In section \ref{proofs} we also present a proof
of a remarkable 1948 result of R. A. Rankin which provides a
beautiful solution to a particular problem dating from 1741. The
first solution to this problem was given in 1886 by W.\ H.\
Thompson. This solution is more subtle than the proofs in section
\ref{leads}.  In section \ref{newthm} we present a small
result concerning $S_4$.

No knowledge of change-ringing is assumed for this
article. We keep the bell-ringing terminology to a minimum,
although a little is necessary. The mathematics involved is
elementary group theory.

\bigskip
\section{The Johnson-Steinhaus-Trotter solution} \label{jst}
\bigskip

The Johnson-Steinhaus-Trotter solution of the motel problem
described  below satisfies rules (1), (2), (3), (4M), but not rule
(4R) for $n>3$, and is therefore not a solution to the ringers
problem for $n>3$.

\bigskip

The solution to the motel problem in Johnson \cite{Johnson},
Steinhaus \cite{Steinhaus} and Trotter \cite{Trotter} (see also
\cite{Gard}, \cite{HallK}, \cite{Lehmer}) is a beautifully simple
idea. We summarise the idea as follows. Construct all arrangements
of $n+1$ objects inductively, using all arrangements of $n$
objects. The induction begins because 1 is a solution for $n=1$,
or $12, 21$ is a solution for $n=2$. Expand the list of all
arrangements of $n$ objects by replacing each arrangement by $n+1$
copies of itself. Place the new $(n+1)$th object $M$ on the left
of the first arrangement. To obtain the next arrangement, we
interchange $M$ with the object to the right of $M$. After doing
this $n$ times, we have $M$ on the extreme right. Then we leave
$M$ here for one step, as the old $n$ objects undergo their first
rearrangement. Next $M$ moves to the left, one place at a time,
and then $M$ stays at the extreme left for one step after arriving
there, as the  old $n$ objects undergo their second rearrangement.
We continue this process until we reach the end of the list.

For example, to obtain the six arrangements of three objects from
the two arrangements of two objects, we first expand and obtain
23, 23, 23, 32, 32, 32. Then we weave 1 through these as
instructed in the algorithm: 123, 213, 231, 321, 312, 132.

To get the solution for four objects, we write the solution for
three objects as 234, 324, 342, 432, 423, 243. Then, with 1 as the
new object $M$, we obtain the list of rows shown below. The reader
should check the list, and ignore the $A$'s, $B$'s and $C$'s for
the moment.

\bigskip

\begin{tabular}{rlllrl}

1\ 2\ 3\ 4 & &\\  & $A$ & & $A$ & & $A$ \\  2\ 1\ 3\ 4 & & 3\ 1\
4\ 2 & & 4\ 1\ 2\ 3\\  & $B$ &          & $B$ & & $B$ \\  2\ 3\ 1\
4 & & 3\ 4\ 1\ 2 & & 4\ 2\ 1\ 3\\  & $C$ & & $C$ & & $C$
\\  2\ 3\ 4\ 1 & & 3\ 4\ 2\ 1 & & 4\ 2\ 3\ 1\\  & $A$ & & $A$
& & $A$ \\  3\ 2\ 4\ 1 & & 4\ 3\ 2\ 1 & & 2\ 4\ 3\ 1\\  & $C$
& & $C$ & & $C$ \\  3\ 2\ 1\ 4 & & 4\ 3\ 1\ 2 & & 2\ 4\ 1\ 3\\
 & $B$ &          & $B$ & & $B$ \\  3\ 1\ 2\ 4 & & 4\ 1\ 3\ 2
& & 2\ 1\ 4\ 3\\  & $A$ & & $A$ & & $A$ \\  1\ 3\ 2\ 4 & & 1\ 4\
3\ 2 & & 1\ 2\ 4\ 3\\  & $C$ & & $C$ & & $C$ \\  1\ 3\ 4\ 2 & & 1\
4\ 2\ 3 & & 1\ 2\ 3\ 4
\end{tabular}

\bigskip

Now we explain the $A, B, C$ notation. The $A$ stands for the
permutation $(1\ 2)$, $B$ for $(2\ 3)$ and  $C$ for $(3\ 4)$,
acting on the positions. Between two rows we have listed the
permutation used to get from one to the other. A shorter way to
write this is simply to list the sequence of $A$'s, $B$'s and
$C$'s,
$$A,B,C,A,C,B,A,C,A,B,C,A,C,B,A,C,A,B,C,A,C,B,A,C,$$
and we will sometimes do this in the sequel.

\bigskip

We may claim that the motel problem was actually solved in the
1600's by bell-ringers. In the early days ringers did not have
rule (4R), and all early methods of changing the order of the
bells involved only two adjacent bells switching place at any one
time. Such changes were called ``plain changes.'' And so, in fact,
rule (4M) was being used instead of rule (4R). Thus, the ringers
problem in those days was identical with the motel problem. An
exact copy of the solution given above to the motel problem for
$n=4$ can be seen in a book dating from c.\ 1621. The full 120
plain changes on five bells were being rung by the mid seventeenth
century.

%\epsfxsize=4in

%\epsfbox{cam3.eps} Fig 1

%\newpage

Fabian Stedman's Campanalogia \cite{Stedman} published in 1677
gives  the Johnson-Steinhaus-Trotter solution of the motel problem
for $n=3, 4, 5$ and $6$. %see figure 2 for $n=4$.

\bigskip

%\epsfxsize=3in

%\epsfbox{cam5.eps}

%\epsfxsize=3in

%\epsfbox{cam10.eps} Fig 2

%\bigskip

Moreover, the general pattern of the algorithm has been noted, as
is evidenced by the following (paraphrased) quote concerning the
$6!$ plain changes on 6 bells.

\bigskip

{\sl The method of the seven hundred and twenty, has an absolute
dependency upon the method of the sixscore changes on five bells;
for five of the notes are to make the sixscore changes, and the
sixth note hunts continually through them, and every time it leads
or lies behind them, one of the sixscore changes must then be
made. The method of the seven hundred and twenty is in effect the
same as that of the sixscore; for as the sixscore comprehended the
twenty-four changes on four, and the six on three; so likewise the
seven hundred and twenty comprehend the sixscore changes on five,
the twenty-four changes on four, and the six changes on three.} --
F.\ Stedman, 1677

\bigskip

Although the general pattern had been observed, plain changes on
more than 6 bells are difficult to ring and bell-ringers
progressed to other ways of ringing. There is an earlier book than
Stedman's, entitled {\it Tinntinnalogia}, written by Richard
Duckworth and published in 1668, which also discusses plain
changes on 3, 4, 5 and 6 bells. This book presumably gives the
Johnson-Steinhaus-Trotter solution to the motel problem, but we
have not been able to confirm this. Nevertheless, the evidence
from Stedman's book would appear to show that their algorithm had
been discovered three hundred years before.  Knuth \cite{Knuth}
states that there is a document written by Mundy dating from 1653
which gives the algorithm.

\bigskip

%The  claim should perhaps be qualified, because ringers are not
%trained mathematicians, and would not write down a solution for
%general $n$.  Only a solution for the number of bells being rung
%would be written down, which would have been typically at most 8.
%The complete $8!+1$ rows on 8 bells were rung for the first
%time in 1761 [F].
%However, general patterns had certainly been observed and
%a solution for general $n$ could have been written down if desired.
%One can find the result for any $n$ stated in [W1] as
%Theorem 3.2, which is a restatement of a result in [R].
%It is stated that the first solution was found in
%1958 by Steinhaus, but in fact the first published solution
%is due to Stedman in 1668/1677.

%This solution has the nice property that
%at each stage a single transposition is applied.

The following observation will be important in the next sections.
The solution to the motel problem given above for $n=4$ is
equivalent to a way of writing down the elements of $S_4$ in an
ordered list, with the property that each element of the list is
obtained from the preceding element by multiplication on the right
by one of $A, B$ or $C$. The list would begin $A, AB, ABC, ABCA,
ABCAC, \ldots $.

%%%%%%%%%%

\bigskip
\section{Plain Hunt} \label{plainhunt}
\bigskip

The simplest method for bell-ringers is called Plain Hunt. We will
describe Plain Hunt on 5 bells, and then the generalisation to $n$
bells. But note that not all methods generalise easily to $n$
bells.

Plain Hunt on 5 bells uses two permutations applied
alternately to rounds, until rounds comes back again.
The permutations are $X=(1\ 2)(3\ 4)$ and
$Y=(2\ 3)(4\ 5)$.  We write $X$ or $Y$ between two
rows to indicate which permutation has been used
to get from one to the other.
After each element $X$ or $Y$ we list the product of all
the elements so far.

\bigskip

\begin{tabular}{rll}
 & & \\
 1\ 2\ 3\ 4\ 5 & & \\
 & $X,\ (1\ 2)(3\ 4)$\\
 2\ 1\ 4\ 3\ 5 & &\\
 & $Y,\ \ XY=(1\ 3\ 5\ 4\ 2)$\\
 2\ 4\ 1\ 5\ 3 & &\\
 & $X,\ \ XYX=(1\ 4)(3\ 5)$\\
 4\ 2\ 5\ 1\ 3 & &\\
 & $Y,\ (1\ 5\ 2\ 3\ 4)$\\
 4\ 5\ 2\ 3\ 1 & &\\
 & $X,(1\ 5)(2\ 4)$\\
 5\ 4\ 3\ 2\ 1 & &\\
 & $Y,(1\ 4\ 3\ 2\ 5)$\\
 5\ 3\ 4\ 1\ 2 & &\\
 & $X,(1\ 3)(2\ 5)$\\
 3\ 5\ 1\ 4\ 2 & &\\
 & $Y,(1\ 2\ 4\ 5\ 3)$\\
 3\ 1\ 5\ 2\ 4 & &\\
 & $X,(2\ 3)(4\ 5)$\\
 1\ 3\ 2\ 5\ 4 & &\\
 & $Y$,identity\\
 1\ 2\ 3\ 4\ 5 & &
\end{tabular}

\bigskip

Using the above shorthand, we could write this as
$$ X, Y, X, Y, X, Y, X, Y, X, Y.$$
After 10 permutations we return to rounds. In other words, there
is a total of 11 rows in Plain Hunt on 5, if we include both
rounds at beginning and end. We could have predicted this; since
$X$ and $Y$ generate a group of order 10 it follows that there can
be no more than 10 permutations in Plain Hunt on 5, and indeed in
any method only using $X$ and $Y$.  However, we should note that
given any two permutations $A$ and $B$, if they generate a group
of order $m$ it does not follow that we can ring a method with
$m+1$ rows using $A$ and $B$. This is because of Rule (3).

In Plain Hunt on 6 bells, the permutations used are
$X=(1\ 2)(3\ 4)(5\ 6)$ and $Y=(2\ 3)(4\ 5)$.
On 7 bells we use
$X=(1\ 2)(3\ 4)(5\ 6)$ and $Y=(2\ 3)(4\ 5)(6\ 7)$.
On 8 bells we use
$X=(1\ 2)(3\ 4)(5\ 6)(7\ 8)$ and $Y=(2\ 3)(4\ 5)(6\ 7)$.
The generalisation to $n$ bells is now clear.

Note above that $X$ and $Y$ are products of disjoint
transpositions of consecutive numbers.   This is
demanded by rule (3).

The permutations $X$ and $Y$ generate a subgroup $H_{2n}$ of $S_n$
of order $2n$, which we will call the {\it hunting subgroup} (even
though it is better known as the dihedral group $D_{2n}$). Note
that if we take a list of the elements after the comma above, we
obtain a list of the elements of the subgroup  $H_{10}$.

\bigskip

{\bf Remark 1}.
Here is the general idea for solving the ringers problem:
to devise a method with more than $2n$ rows,
we must throw another permutation into the mix.
We wish to use as few permutations as possible
in order to keep the method as simple as possible,
while obtaining a method as long as possible,
hopefully with all $n!+1$ rows.

The first solutions employed by ringers to the ringers problem
involved cosets of the hunting/dihedral subgroup $H_{2n}$ of $S_n$
generated by $X=(1\ 2)(3\ 4)\ldots (n-1,\ n)$ and $Y=(2\ 3)(4\
5)\ldots (n-2,\ n-1)$ when $n$ is even, and $X=(1\ 2)(3\ 4)\ldots
(n-2,\ n-1), Y=(2\ 3)(4\ 5)\ldots (n-1,\ n)$ when $n$ is odd.
These particular involutions are used for the transitions between
successive rows because rules (3) and (4R) will be obeyed.
%As far as ringers are concerned,
%a transition from one row to the next will always be
%a product of disjoint transpositions of consecutive numbers.
As long as we do not apply $X$ or $Y$ twice in succession, rule
(2) will be obeyed.

\bigskip

{\bf Remark 2}. Let us mention here another rule, which roughly
states that each bell follows the same path. We will not go into
any further detail on this. One can see that this is indeed the
case in Plain Hunt, and that it will also hold in a coset of the
hunting subgroup. This is why we shall assume for this article
that all methods are a union of cosets of the hunting subgroup.
Using cosets keeps the method simple, one of the goals from Remark
1. How to choose the cosets such as to obey rules (1)-(3) and (4R)
is the real question.

\bigskip
\section{Plain Bob}\label{plainbob}
\bigskip

Probably the next simplest method after Plain Hunt is called
Plain Bob.
This method dates from about 1650.
We will describe Plain Bob on 4 bells,
and then 6 bells.
The idea of Plain Bob is to combine Plain Hunt with some
particular cosets.

We do this because Plain Hunt is not yet a solution to the ringers
problem of finding an extent of all $n!+1$ rows.  We obtain a
solution to the ringers problem for $n=4$ using the hunting
subgroup $H_8$, a group of order 8, and two left cosets, which are
$(2\ 4\ 3)H_8$ and $(2\ 3\ 4)H_8$. The advantage of using cosets
of a subgroup is that distinct cosets are disjoint and therefore
rule (2) is automatically satisfied. This fact was surely known
to, and utilised by, the early composers. Here is the full
solution, which uses one other permutation, namely $Z=(3\ 4)$,
apart from $X$ and $Y$.  This $Z$ is used in order to switch into
the cosets.

We spell this out in detail.  First consider Plain Hunt on 4
bells:
\bigskip

\begin{tabular}{rll}

 & & \\  1\ 2\ 3\ 4 & & \\  & $X,(1\ 2)(3\ 4)$\\  2\ 1\
4\ 3 & &\\  & $Y,(1\ 3\ 4\ 2)$\\  2\ 4\ 1\ 3 & &\\  & $X,(1\ 4)$\\
4\ 2\ 3\ 1 & &\\  & $Y,(1\ 4)(2\ 3)$\\  4\ 3\ 2\ 1 & &\\  & $X,(1\
3)(2\ 4)$\\  3\ 4\ 1\ 2 & &\\  & $Y,(1\ 2\ 4\ 3)$\\  3\ 1\ 4\ 2 &
&\\  & $X,(2\ 3)$\\  1\ 3\ 2\ 4 & &\\  & $Y$, identity\\  1\ 2\ 3\
4 & &
\end{tabular}

\bigskip

where $X=(1\ 2)(3\ 4)$ and $Y=(2\ 3)$. The sequence of elements
after the commas consists of the elements of $H_8$.

In Plain Bob on 4 bells, instead of doing the final $Y$ which
takes us back to rounds, we do $Z=(3\ 4)$ instead. We then
continue with $X$ and $Y$ alternately. This has the effect of
taking us into the left coset $(Y^{-1}Z)H_8$ of the hunting
subgroup $H_8$, i.e., after the comma we will be listing the
elements of $(Y^{-1}Z)H_8$. Here is what happens.

\begin{tabular}{rll}
 & & \\
 3\ 1\ 4\ 2 & &\\
 & $X,\ XYXYXYX=(2\ 3)= Y^{-1}$\\
 1\ 3\ 2\ 4 & &\\
 & $Z,\ Y^{-1}Z=(2\ 3)(3\ 4)=(2\ 4\ 3)$\\
 1\ 3\ 4\ 2 & &\\
 & $X,\ Y^{-1}ZX$=(1\ 2\ 3) \\
 3\ 1\ 2\ 4 & &\\
 & $Y,\ Y^{-1}ZXY$=(1\ 3) \\
 3\ 2\ 1\ 4 & &\\
  \hfil\vdots  &  \vdots  & \hfill \vdots

\end{tabular}

\bigskip

At the same point in this coset when we have reached
$Y^{-1}ZXYXYXYX=$ \break $Y^{-1}ZY^{-1}$, instead of doing the
final $Y$ next (which would cause us to repeat $1\ 3\ 4\ 2$,
disobeying rule (2)) we do $Z$ again, which takes us into the
coset $(Y^{-1}Z)^2H_8$. Again we alternate between $X$ and $Y$ to
take us through this coset, and at the same point again we do $Z$,
which takes us back to rounds. Here is the full set of rows.

\bigskip

\begin{tabular}{rlllrl}
 1\ 2\ 3\ 4 & &\\
 & $X,(1\ 2)(3\ 4)$ &          & $X,(1\ 2\ 3)$ & & $X,(1\ 2\ 4)$ \\
 2\ 1\ 4\ 3 & & 3\ 1\ 2\ 4 & & 4\ 1\ 3\ 2\\
 & $Y,(1\ 3\ 4\ 2)$ &          & $Y,(1\ 3)$ & & $Y,(1\ 3\ 2\ 4)$ \\
 2\ 4\ 1\ 3 & & 3\ 2\ 1\ 4 & & 4\ 3\ 1\ 2 \\
 & $X,(1\ 4)$ &                & $X,(1\ 4\ 3\ 2)$ & & $X,(1\ 4\ 2\ 3)$ \\
 4\ 2\ 3\ 1 & & 2\ 3\ 4\ 1 & & 3\ 4\ 2\ 1\\
 & $Y,(1\ 4)(2\ 3)$ &          & $Y,(1\ 4\ 2)$ & & $Y,(1\ 4\ 3)$ \\
 4\ 3\ 2\ 1 & & 2\ 4\ 3\ 1 & & 3\ 2\ 4\ 1\\
 & $X,(1\ 3)(2\ 4)$ &          & $X,(1\ 3\ 4)$ & & $X,(1\ 3\ 2)$ \\
 3\ 4\ 1\ 2 & & 4\ 2\ 1\ 3 & & 2\ 3\ 1\ 4\\
 & $Y,(1\ 2\ 4\ 3)$ &          & $Y,(1\ 2\ 3\ 4)$ & & $Y,(1\ 2)$ \\
 3\ 1\ 4\ 2 & & 4\ 1\ 2\ 3 & & 2\ 1\ 3\ 4\\
 & $X,(2\ 3)$ &                & $X,(2\ 4)$ & & $X,(3\ 4)$ \\
 1\ 3\ 2\ 4 & & 1\ 4\ 3\ 2 & & 1\ 2\ 4\ 3\\
 & $Z,(2\ 4\ 3)$ &             & $Z,(2\ 3\ 4)$ & & $Z$,identity \\
 1\ 3\ 4\ 2 & & 1\ 4\ 2\ 3 & & 1\ 2\ 3\ 4
 \end{tabular}

\bigskip

This solution could be represented by the sequence of permutations
$$X,Y,X,Y,X,Y,X,Z,X,Y,X,Y,X,Y,X,Z,X,Y,X,Y,X,Y,X,Z.$$

There are a few important observations to make here. Firstly, by
adding cosets of $H_8$ to Plain Hunt, we have increased the number
of rows in our set of changes. This is the basic idea of all the
methods considered in this paper, as we said in Remark 2.

Secondly, we have listed (after the comma)
the permutations in the order
$$H_8\backslash \{\hbox{identity}\}, \ (Y^{-1}Z)H_8, \ (Y^{-1}Z)^2H_8,
\ \hbox{identity}$$
where within each coset we use the order from Plain Hunt.

Thirdly, we note that we obtained 3 cosets in total
because $Y^{-1}Z=(2\ 4\ 3)$ has order 3.
We also note that the union of the 3 cosets is all
of $S_4$, so in this case we obtained the maximum
number of permutations possible, an extent.
This does not happen in general.

\bigskip

Next we summarise Plain Bob on 6 bells.   Recall that the idea is
to add cosets of $H_{12}$ to Plain Hunt on 6, so that we obtain
more rows.  Of course, this must be done without disobeying any of
rules (1), (2), (3), (4R). We recall that the generators of the
hunting subgroup are $X=(1\ 2)(3\ 4)(5\ 6)$ and $Y=(2\ 3)(4\ 5)$.
In addition we use $Z=(3\ 4)(5\ 6)$, in the same way as in Plain
Bob on 4 bells. Since $Y^{-1}Z=(2\ 3)(4\ 5)(3\ 4)(5\ 6)=(2\ 4\ 6\
5\ 3)$ has order 5, we obtain 5 cosets of $H_{12}$ for a total of
$60$ permutations (61 rows including both rounds).

Let us check that rules (1)--(3) and (4R) are satisfied. Because
$X, Y, Z$ are all products of disjoint transpositions of
consecutive numbers and they are the only permutations used to get
from one row to the next, rule (3) is satisfied. Rule (1) is also
satisfied because of the construction, and rule (2) is satisfied
because the cosets are distinct, and therefore disjoint. Rule (4R)
is satisfied because it is satisfied for Plain Hunt, and Plain Bob
is a union of cosets of Plain Hunt. Here we see that group theory
provides us with a construction of a method, and a proof that it
obeys the rules, without writing out all the rows.

\bigskip

One can also ring Plain Bob on an odd number of bells. For any
$n$, Plain Bob on $n$ bells uses $n-1$ cosets of $H_{2n}$ and so
has $2n(n-1)$ permutations. Only when $n=4$ does $2n(n-1)=n!$.

This solution to the ringers problem dates from the 17th century,
and can be found in Stedman's 1677 book \cite{Stedman}. %see figure
%3 below.

%\epsfxsize=3in

%\epsfbox{cam8.eps} Fig 3

%Since any solution to the ringers problem is also a solution to
%the motel problem, we have
%a solution to the motel problem in the $n=4$ case.
%Note that rules (1)-(4)
%are obeyed, so this is a different solution to the one in section 3.
Note that bell 1 here behaves in the same way as in the section
\ref{jst}  solution to the motel problem. The other bells behave
differently however. Unfortunately, this solution to the ringers
problem does not generalise to arbitrary $n$ in an obvious way,
unlike the motel problem. The task of combining cosets of $H_{2n}$
to obtain a solution to the ringers problem for $n>4$ is highly
nontrivial, and solutions involve many clever ideas. See
\cite{WhiteRing} for more details.

\section{A New Method}\label{ccc}

We now introduce a "new" method on 5 bells which is a union of
cosets of $H_{10}$.  We mimic the construction of Plain Bob on 5
bells except that instead of using $Z=(3\ 4)$ we will use $Z= (1\
2)$.  This results in 6 cosets of $H_{10}$ because $Y^{-1}Z=(2\
3)(4\ 5)(1\ 2)=(1\ 2\ 3)(4\ 5)$ has order 6.  Thus we obtain a
method with 61 rows (including both rounds).

This method is not listed in the collection of known methods,
so we propose to call it Christ Church
Dublin Differential Doubles. 
This is now recognised, see the web site \cite{Smith}
under differentials.
We outline a bob for this method in
section \ref{leads}.

\bigskip
\section{A Group-Theoretic Formulation}\label{groupform}
\bigskip

As explained at the end of section \ref{jst}, the solution to the
ringers problem given above is equivalent to a way of writing down
the elements of $S_4$ in a list, with the property that each
element of the list is obtained from the preceding element by
multiplying by one of $X, Y$ or $Z$. This idea motivates the
following definition.

\bigskip
\subsection{Unicursal Generation}
\bigskip

{\bf Definition}. Let $G$ be a finite group of order $n$, and let
$T$ be a subset of $G$. We say that $T$ generates $G$
$\underline{unicursally}$ if the elements of $G$ can be ordered
$g_1, g_2, \ldots ,g_n$ so that for each integer $i$, there exists
$t_i\in T$ such that $g_{i+1}=g_it_i$. (Here subscripts are
considered modulo $n$.)

In this framework, the motel problem can be restated as follows:
is $S_n$ generated unicursally by $T=\{(1\ 2), (2\ 3), (3\ 4),
\ldots ,(n-1,\ n)\}$ ?

The ringers problem can be restated as follows: is $S_n$ generated
unicursally by a subset $T$ satisfying the following conditions:
\begin{enumerate}
\item each element of $T$ is a product of disjoint transpositions
of consecutive numbers (this is rule 3)  \item $t_i$ and $t_{i+1}$
have no common fixed point for any $i$ (this is rule 4R).
\end{enumerate}

\bigskip

{\bf Example 1.}  Let $G=S_3$, and let $T=\{(1\ 2),(2\ 3)\}$. Then
$T$  generates $G$ unicursally. The reader will find it useful to
verify this small example.
\bigskip

{\bf Example 2.}  Let $G$ be the $D_4$ from section
\ref{plainhunt}, and let $T=\{(1\ 2)(3\ 4), (2\ 3)\}$. Then $T$
generates $G$ unicursally as shown in section \ref{plainhunt}.
\bigskip

{\bf Example 3.}  Let $G=S_4$, and let $T=\{(1\ 2)(3\ 4), (2\ 3),
(3\ 4)\}$ as in section \ref{plainbob} in the solution to the
ringers problem for $n=4$.
\bigskip

{\bf Example 4.}  Let $G=S_n$.  The solution to the the motel
problem described in section \ref{jst} shows that $S_n$ is
generated unicursally by the $n-1$ elements of $T=\{(1\ 2), (2\
3), (3\ 4), \ldots ,$ $(n-1,\ n)\}$.
\bigskip

{\bf Example 5.} Let $G$ be any finite group and let $T=G$. Then
$G$ is generated unicursally by $T$, and any ordering of the
elements of $G$ can be used, since $g(g^{-1}h)=h$ for any $g, h\in
G$.
\bigskip

{\bf Example 6.} Let $G=S_n$ and let $$T=\{(1\ 2), (1\ 2)(3\ 4)(5\
6)\ldots , (2\ 3)(4\ 5)(6\ 7)\ldots \}.$$ It is shown in
\cite{Rapaport} (see also \cite{WhiteRing}) that $T$ generates $G$
unicursally.  Knuth \cite{Knuth} states that Rapaport's result has
been generalised by Savage.

\bigskip

%As the subscripts are modulo $n$ in the definition,
%without loss of generality we may assume $g_1=1$.

{\bf Remarks}. In order for a subset $T$ to generate $G$
unicursally it is necessary that $T$ generates $G$. This condition
is not sufficient, as examples 7, 8 and 9 below show.

We are usually interested in the case when $T$ is a small, and
often minimal, set of generators. The general question of  whether
a given $G$ is generated unicursally by a given $T$ seems very
difficult. This problem may be related to word problems in the
group $G$.

Given a generating set $T$ for a finite group $G$, the Cayley
colour graph $C_T(G)$ is the graph with the elements of $G$ for
vertices, and all directed edges $(x,xt)$ where $t\in T$. Each
directed edge is coloured by the generator $t$.  If every element
of $T$ has order 2, then the graph may be considered undirected.
Usually assumptions are made to ensure that $C_T(G)$ has no loops
or multiple edges.  The group $G$ acts regularly and transitively
on the vertices of $C_T(G)$, and is the automorphism group of
$C_T(G)$.  The following theorem is clear.

\begin{theorem}
A group $G$ is generated unicursally by $T$ if and only if the
Cayley colour graph defined by $G$ and $T$ is Hamiltonian.
\end{theorem}

This is the point of view of Rapaport \cite{Rapaport} and White
\cite{WhiteCamb}, \cite{WhiteRing}, \cite{WhiteSted}. White has
written several papers on bells and topological graph theory.

\bigskip

The motel problem and the ringers problem are concerned with
specific types of subset $T$ of the symmetric group $S_n$. For the
motel problem, as we have said above, $T$ is the set $\{(1\ 2),
(2\ 3), (3\ 4), \ldots (n-1,\ n)\}$. For the ringers problem,
elements of $T$ will be products of disjoint transpositions of
consecutive numbers (because of rule 3), and one must ensure that
$t_i$ and $t_{i+1}$ have no common fixed point (because of rule
4R).  Which particular $T$ is chosen depends on the method.
Two methods will be discussed in this paper, Plain Bob
and Grandsire.

Let us generalise and consider the unicursal generation of
$S_n$ by elements other than products of disjoint transpositions.
First we make a few simple observations about arbitrary groups.
The classification of groups generated by $T$ of size 1 is
straightforward.

\begin{theorem}
A group $G$ is generated unicursally by a subset
of size 1 if and only if $G$ is cyclic.
\end{theorem}

Proof. Suppose $G$ is generated unicursally by $T=\{x\}$. Assuming
$g_1=1$ (w.l.o.g.) then $g_2=x$, and then $g_3=x^2$, and
$g_{i+1}=g_ix=x^i$ for any $i$. The other implication is clear.

\bigskip

%A chain of length $r$ with respect to $T$
%is any sequence of elements
%$g_1, g_2, \ldots ,g_r$ so that for each $i$,
%there exists $t_i\in T$ such that
%$g_{i+1}=g_it_i$
%(subscripts are taken modulo $r$).

The classification of groups generated by $T$ of size 2 is
nontrivial. Clearly if $G$ is isomorphic to a direct product of
two cyclic groups then $G$ is unicursally generated by a subset of
size 2. The following theorem is a remarkable result of R.\ A.\
Rankin on this case.
%It allows us to easily construct examples of groups $G$
%and subsets $T$ of size 2 with $G$ not
%generated unicursally by $T$.
The result in Rankin's paper \cite{Rankin} is more general than
the version we state here, and is somewhat based on ideas of
Thompson \cite{Thompson}. The end result is very simply stated in
group theoretic language, even though the problem and the proof
are somewhat combinatorial. We give a proof in section
\ref{proofs}.

Let $\langle g \rangle$ denote the subgroup generated by $g$.

\begin{theorem}  \label{rankthm}
(Rankin, 1948) Let $G$ be a finite group. Suppose that $G$
is generated by $T=\{x,y\}$, and that $\langle x^{-1}y \rangle$
has odd order. If $G$ is generated unicursally by $T$, then
$\langle x \rangle$ and $\langle y \rangle$ have odd index in $G$.
\end{theorem}

\bigskip

{\bf Example 7.}  It is easily checked that $G=A_4$ is generated
by $T=\{A, B\}$ where $A=(1\ 2\ 3), B=(1\ 2\ 4)$. However, $A_4$
is not generated unicursally by $T$ by Rankin's theorem, because
$A^{-1}B=(1\ 3\ 4)$ has odd order but $\langle A \rangle$ has
index 4 in $A_4$.
%The longest chain generated unicursally
%by $T$ is of length 11.
%One such chain is
\bigskip

{\bf Example 8.} It is easily checked that $A_5$ is generated by
$T=\{A, B\}$, where $A=(1\ 3\ 5\ 4\ 2)$ and $B=(3\ 5\ 4)$. Here
$A^{-1}B$ has order 3, so Rankin's theorem applies. Since $\langle
B \rangle$ has even index in $A_5$ we conclude that $A_5$ is not
generated unicursally by $T$.

\bigskip

{\bf Example 9.} It is not hard to show that $S_n$ is generated by
the transposition $A=(n-1,\ n)$ and the $(n-1)$-cycle $B=(1\ 2\ 3\
\ldots n-1)$. We may well ask whether $S_n$ is  generated
unicursally by $A$ and $B$. Suppose $n$ is odd and $n>3$.   Then
$S_n$ is not generated unicursally by $T=\{A, B\}$ by Rankin's
theorem, because $A^{-1}B=(1\ 2\ \ldots n)$ has odd order but
$\langle A \rangle$ has even index in $S_n$ for $n>3$.

\bigskip

{\bf Example 10.} It is not hard to show that $S_n$ is generated
by $\sigma=(1\ 2\ \ldots n)$ and $\tau = (1\ 2)$.
We (of course) ask if $S_n$ is generated unicursally by these
elements.
Rankin's theorem shows that the answer is negative if $n\geq 4$ is
even, since $\tau^{-1} \sigma$ is an $(n-1)$-cycle.
We will mention this example again soon.

\bigskip

From the discussion in the previous sections, the following is now
obvious (and has been observed before, see \cite{HallK} for
example). Let $t_i$ be as above.

\begin{theorem}
\begin{enumerate}
\item  The existence of an extent on $n$ bells
satisfying rules (1)-(3), where the allowed permutations between
rows are $X_1, \ldots ,X_k$, is equivalent to $T=\{X_1, \ldots
,X_k\}$ generating $S_n$ unicursally. \item  The existence of an
extent on $n$ bells satisfying rules (1),(2),(3),(4R), where the
allowed permutations between rows are $X_1, \ldots ,X_k$, is
equivalent to $T=\{X_1, \ldots ,X_k\}$ generating $S_n$
unicursally with the additional property that $t_i$ and $t_{i+1}$
have no common fixed point for any $i$.
\end{enumerate}
\end{theorem}

A permutation $X_i$ being an allowed transition between rows is
equivalent to $X_i$ being a product of disjoint transpositions of
consecutive numbers, by rule (3). No bell staying in the same
place for more than two rows (rule (4R)) is equivalent to no two
consecutive transitions having a common fixed point.

%\proclaim{Corollary 2}.
%If $T=\{X_1, \ldots ,X_k\}$ does not generate $S_n$
%unicursally, then there does not exist an extent
%on $n$ bells using $X_1, \ldots ,X_k$ as permutations
%between rows.

As we mentioned in Remark 2 at the end of section 3, we assume
that the methods in this article are cosets of $H_{2n}$.
Therefore, included in $T$ will be $X$ and $Y$, the generators for
the hunting subgroup $H_{2n}$. But note that even if $S_n$ is
generated unicursally by $T$, it does not follow that it is
generated unicursally as cosets of $H_{2n}$.

\bigskip

{\bf Remark 3}. In the ringing methods we discuss in this paper,
one can divide up an extent on $n$ bells into groups of $2n$ rows
called {\it leads}, roughly (but not exactly) corresponding to
cosets of $H_{2n}$. A method composed of cosets of $H_{2n}$ is
also a method made up of a succession of leads. In this paper we
use only two types of leads, and we consider methods and extents
made from a sequence of these leads. In these cases then, one can
show that the existence of an extent on $n$ bells of this type is
equivalent to the alternating group $A_{n-1}$ being generated
unicursally by $T'$, where $T'$ is a set of generators related to
$T$. For more details on this, see section 6.

The general question of  whether
a given $G$ is generated unicursally by a given $T$
seems very difficult.

\bigskip
\subsection{The famous old question}
\bigskip

The famous old question mentioned in the introduction concerns the
three permutations $X=(1\ 2)(3\ 4)(5\ 6), Y=(2\ 3)(4\ 5)(6\ 7)$
and $Z=(1\ 2)(4\ 5)(6\ 7)$ in $S_7$, and whether $S_7$ is
generated unicursally by $X, Y, Z$ in a particular way.  We will
explain this in detail in section \ref{leads}.

This was asked in 1741 by a bell-ringer John Holt, who was able to
construct a method of $4998$ permutations, but could not obtain a
method of $7!=5040$ permutations.  He then (naturally!) queried
the existence of such an extent. As in Remark 3, it can be shown
(see \cite{Rankin},  or section \ref{leads}) that this question is
equivalent to:

\bigskip

{\bf Question A}. Is $A_6$ (acting on $\{2,3,4,5,6,7\}$) generated
unicursally by the two permutations $(3\ 4\ 6\ 7\ 5)$ and $(2\ 4\
7)(3\ 6\ 5)$ ?

\bigskip

The first proof that the answer is no is due to Thompson 
\footnote{Thompson was a civil servant in India at the time,
and a Cambridge mathematics gradute.} (1886)
\cite{Thompson}, with some case-by-case analysis. An insightful
proof was given by Rankin (1948) \cite{Rankin}, where he came up
with theorem \ref{rankthm} based somewhat on Thompson's ideas. 
See also \cite{Rankin2} and a proof by Swan \cite{Swan}.
We present Rankin's proof in section \ref{proofs}, in our special
case only. Rankin's result is more general. Most of the ideas of
the proof can be found in our proof of the special case in section
\ref{proofs}. We also show in section \ref{proofs} that $4998$ is
best possible.

\bigskip

\subsection{Open Questions}

\bigskip

The first concerns example 9.  The argument there works when $n$
is odd, so it is natural to inquire as to what happens when $n$ is
even. Thus, we wonder whether $S_n$ is generated unicursally by
$A=(n-1,\ n)$ and $B=(1\ 2\ 3\ \ldots n-1)$. Rankin's theorem does
not apply directly. However, in the $n=4$ case, by modifying the
argument in the proof of Rankin's theorem, we will show (see
section \ref{newthm}) that $S_4$ is not generated unicursally by
$A=(1\ 2\ 3)$ and $B=(3\ 4)$. The question for even $n\geq 6$
remains open, as far as we are aware.

\bigskip

{\bf Problem 1}:  Let $n\geq 6$ be even. Is $S_n$  generated
unicursally by $T=\{A,B\}$ where $A=(n-1,\ n)$ and $B=(1\ 2\ 3\
\ldots n-1)$ ?

\bigskip
According to Knuth \cite{Knuth} a similar question
was asked in 1975 by Nijenhuis and Wilf in their book
\emph{Combinatorial Algorithms}.  They asked if $S_n$ is generated
unicursally by $\sigma=(1\ 2\ \ldots n)$ and $\tau = (1\ 2)$.
Rankin's theorem shows that the answer is negative if $n\geq 4$ is
even (see example 10).  
Recently, Ruskey-Jiang-Weston \cite{rjw} did a computer
search for $n=5$ and did find that $S_5$ IS unicursally generated
by $\sigma$ and $\tau$. Thus this question is different to problem
1.

\bigskip

This example is not relevant to bells  because the
generators are not products of disjoint transpositions.
However, as the answer to problem 1 may well be negative,
the discussion raises the natural question of whether
$S_n$ is generated unicursally by any two of its elements.
(Example 4 shows that $S_n$ is generated unicursally by
three of its elements.)
We have found that $S_4$ IS generated unicursally by
$A=(1\ 2\ 3)$ and $B=(1\ 2\ 3\ 4)$.
Here is one listing which does the trick:
$$B,A,B,A,A,B,B,B,A,B,B,A,A,B,B,B,A,A,B,A,B,B,B,A.$$
The question for  $n\geq 5$ remains open.

\bigskip

{\bf Problem 2:}  Is $S_n$  generated
unicursally by some two elements for all $n$?

\bigskip

By the above comments, the $n=5$ case has been done.

{\bf Remarks}. A similar eighteenth-century problem to our famous
old question, concerning another method called Stedman, remained
unsolved until 1995. We may discuss this in a future article.
Right transversals of $PSL(2,5)$ in $S_7$ are used to construct
extents in this problem, see \cite{WhiteRing}.

Readers interested only in the proof of the answer to question A
should skip ahead to section \ref{proofs}.

\bigskip
\section{Leads}\label{leads}
\bigskip

In this section we shall explain in detail the two types of leads
mentioned in Remark 3. We only deal in this section with the
methods of Plain Bob on 4 and 6 bells, and Grandsire on 5 and 7
bells, although our remarks have wider application. Then we shall
explain the origin of the famous old question.

\subsection{Plain Bob}

Consider the 25 rows in Plain Bob on 4 bells in section 4; the
first one is rounds, and then the rows can be divided into three
sets of eight. Each of these sets of eight is called a {\it lead}.
In each of these leads, note that each bell is twice in the first
position. Also note that bell number 1 is always in the first
position in the last two rows of each lead. The second of these
two rows, which is the last row of the lead, is called the {\it
lead head}. This holds in general, for Plain Bob on $n$ bells,
where leads have $2n$ rows. In Plain Bob on 4 bells, the lead
heads are $1\ 3\ 4\ 2$, $1\ 4\ 2\ 3$, and $1\ 2\ 3\ 4$. The
following are simple observations:

\begin{enumerate}
\item  The first, second and third lead heads are the result of
$P=Y^{-1}Z=(2\ 4\ 3)$, $P^2$ and $P^3$ respectively acting on
rounds. \item  When considering only lead heads, we may drop the 1
in the first position. \item Plain Bob  can  then be described by
elements of $S_{n-1}$ acting on lead heads.
\end{enumerate}

\bigskip

Now we can fully describe Plain Bob on 6 bells by its lead heads.
Here is the first lead (with initial rounds included as well):

\begin{tabular}{rll}
 & & \\
 1\ 2\ 3\ 4\ 5\ 6& & \\
 & $X, (1\ 2)(3\ 4)(5\ 6)$\\
 2\ 1\ 4\ 3\ 6\ 5 & &\\
 & $Y,(1\ 3\ 5\ 6\ 4\ 2)$\\
 2\ 4\ 1\ 6\ 3\ 5 & &\\
 & $X,(1\ 4)(3\ 6)$\\
 4\ 2\ 6\ 1\ 5\ 3 & &\\
 & $Y, (1\ 5\ 4)(2\ 3\ 6)$\\
 4\ 6\ 2\ 5\ 1\ 3 & &\\
 & $X,(1\ 6)(2\ 4)(3\ 5)$\\
 6\ 4\ 5\ 2\ 3\ 1 & &\\
 & $Y,(1\ 6)(2\ 5)(3\ 4)$\\
 6\ 5\ 4\ 3\ 2\ 1 & &\\
 & $X,(1\ 5)(2\ 6)(3\ 4)$\\
 5\ 6\ 3\ 4\ 1\ 2 & &\\
 & $Y,(1\ 4\ 5)(2\ 6\ 3)$\\
 5\ 3\ 6\ 1\ 4\ 2 & &\\
 & $X,(1\ 3)(2\ 5)(4\ 6)$\\
 3\ 5\ 1\ 6\ 2\ 4 & &\\
 & $Y,(1\ 2\ 4\ 6\ 5\ 3)$\\
 3\ 1\ 5\ 2\ 6\ 4 & &\\
 & $X,(2\ 3)(4\ 5)$\\
 1\ 3\ 2\ 5\ 4\ 6 & &\\
 & $Z, (2\ 4\ 6\ 5\ 3)$\\
 1\ 3\ 5\ 2\ 6\ 4 & &
\end{tabular}

\bigskip

In this case $P=Y^{-1}Z=(2\ 4\ 6\ 5\ 3)$ and the lead heads are
$3\ 5\ 2\ 6\ 4$, $5\ 6\ 3\ 4\ 2$, $6\ 4\ 5\ 2\ 3$, $4\ 2\ 6\ 3\ 5$
and $2\ 3\ 4\ 5\ 6$, corresponding to $P, P^2, P^3, P^4$ and $P^5$
respectively, acting on $2\ 3\ 4\ 5\ 6$. Each of these leads is
called a {\it plain} lead. Each lead head is obtained from the
previous lead head by applying $P$. There are five leads because
$P$ has order 5. This sequence of five plain leads is called a
{\it plain course}.

\bigskip

As we said in Remark 3, there are only two types
of lead considered in this paper.
Let us now describe the other type of lead, at least as
far as Plain Bob on 6 bells goes.
Mathematically there is no reason to have a method
made of only two or three
types of lead, but this is usually what is done
in practice for simplicity and historical reasons
(recall Remark 1).
The complete method is made up of a succession of
leads.

Any plain lead may be described by the sequence of permutations
$$X, Y, X, Y, X, Y, X, Y, X, Y, X, Z.$$
Alternatively, considering only lead heads, we describe a plain
lead by $P$, and the plain course (which has 60 permutations) by
the sequence of five plain leads
$$P, P, P, P, P.$$

The other type of lead is called a {\it bob} lead and
may be described by the
sequence of permutations
$$X, Y, X, Y, X, Y, X, Y, X, Y, X, W$$
where $W=(2\ 3)(5\ 6)$.
If this were done from rounds we would get

\begin{tabular}{rll}
 & & \\
 1\ 2\ 3\ 4\ 5\ 6& & \\
 & $X,\ (1\ 2)(3\ 4)(5\ 6)$\\
 2\ 1\ 4\ 3\ 6\ 5 & &\\
 & $Y,(1\ 3\ 5\ 6\ 4\ 2)$\\
 2\ 4\ 1\ 6\ 3\ 5 & &\\
 & $X,(1\ 4)(3\ 6)$\\
 4\ 2\ 6\ 1\ 5\ 3 & &\\
% & $Y, (1\ 5\ 4)(2\ 3\ 6)$\\
% 4\ 6\ 2\ 5\ 1\ 3 & &\\
% & $X,(1\ 6)(2\ 4)(3\ 5)$\\
% 6\ 4\ 5\ 2\ 3\ 1 & &\\
% & $Y,(1\ 6)(2\ 5)(3\ 4)$\\
% 6\ 5\ 4\ 3\ 2\ 1 & &\\
% & $X,(1\ 5)(2\ 6)(3\ 4)$\\
% 5\ 6\ 3\ 4\ 1\ 2 & &\\
% & $Y,(1\ 4\ 5)(2\ 6\ 3)$\\
% 5\ 3\ 6\ 1\ 4\ 2 & &\\
% & $X,(1\ 3)(2\ 5)(4\ 6)$\\
  \ \ \ \ \vdots & \vdots & \\
 & & \\
 3\ 5\ 1\ 6\ 2\ 4 & &\\
 & $Y,(1\ 2\ 4\ 6\ 5\ 3)$\\
 3\ 1\ 5\ 2\ 6\ 4 & &\\
 & $X,(2\ 3)(4\ 5)$\\
 1\ 3\ 2\ 5\ 4\ 6 & &\\
 & $W,(4\ 6\ 5)$\\
 1\ 2\ 3\ 5\ 6\ 4 & &
 \end{tabular}

\bigskip

This is a bob lead.

%This has the effect of putting us in the coset
%$BH_6$ where $B=Y^{-1}W=(2\ 3)(4\ 5)(2\ 3)(5\ 6)=(4\ 6\ 5)$.

As we used $P$ to denote a plain lead we shall use
$B$ to denote a bob lead.
We can now construct longer methods using
a combination of plain and bob leads.
Here is one such method:
$$P, P, P, P, B, P, P, P, P, B, P, P, P, P, B.$$

This corresponds to doing the first 59 of the 60 permutations in
the plain course. The 60th permutation that would be performed in
a plain course is $Z$, which would bring us back to rounds.
Instead of this last $Z$, we do $W=(2\ 3)(5\ 6)$. This has the
effect of putting us into another coset, namely $BH_{12}$ where
$B=Z^{-1}W=(3\ 4)(5\ 6)(2\ 3)(5\ 6)=(2\ 3\ 4)$. We then repeat the
same 59 permutations, then do $W$ again, then the 59 and then $W$
again, which returns us to rounds since $B$ has order 3. We finish
up with a method of 180 permutations. By a similar argument as in
section \ref{plainbob}, rules (1),(2),(3),(4R), are obeyed.
%The permutation $W$ is called a `bob' (see section 7 for why).

\bigskip

We have still not succeeded in getting an extent
of $6!=720$ permutations.
It is possible that some other
sequence of plain and bob leads will give us an extent.
The following result ends all hope of this.
%More general results can be found in [7].

\begin{theorem} \label{pbminor}
There does not exist an extent of Plain Bob on 6 bells using
plain and bob leads. The longest possible method using plain and
bob leads has $360$ permutations, and there does exist such a
method.
\end{theorem}

Proof:  The key to the proof  is to observe that $P, B, Z$ and $W$
are all even permutations.  The fact that $P$ and $B$ are even
implies that any lead head will be an even permutation of $2\ 3\
4\ 5\ 6$. Also, the row before a lead head is the result of
applying either $Z^{-1}$ or $W^{-1}$ to the lead head. Since $Z$
and $W$ are even, we see that in any method of plain and bob leads
all rows with bell 1 in the first position are followed by an even
permutation of $2\ 3\ 4\ 5\ 6$. The result follows, because if we
did obtain an extent we would get all possible permutations of $2\
3\ 4\ 5\ 6$ following 1.

This argument also shows that any method using plain
and bob leads has at most $5!/2=30$ leads, since each lead
has two rows with 1 in the first position, and these
rows must be followed by an even
permutation of $2\ 3\ 4\ 5\ 6$.
Each lead has 12 rows, so a method
with plain and bob leads has at most
$12\times 30 = 360$ permutations.

To show that $360$ is possible we give an ordering:
$$B, P, P, P, B, B, P, P, P, P,
B, P, P, P, B, B, P, P, P, P, B, P, P, P, B, B, P, P, P, P.$$ The
reader may check that this sequence of plain and bob leads obeys
rules (1),(2),(3),(4R).

\bigskip

{\bf Remarks}.
The use of plain and bob leads applies to Plain Bob
on $n$ bells (and other methods).
The number of leads in an extent on $n$ bells is
$n!/(2n)=(n-1)!/2$, which is the cardinality of $A_{n-1}$.
If $P$ and $B$ generate $A_{n-1}$ unicursally, then
we can construct an extent made up of plain and bob leads.

On 6 bells it is true that $P$ and $B$ generate $A_5$, but example
4 and theorem \ref{pbminor} both show that they do not generate
$A_5$ unicursally. Example 4 used Rankin's theorem, but theorem
\ref{pbminor} gives a different and shorter proof. The argument in
theorem \ref{pbminor} is shorter because parity can be used to
answer the question.  The famous old question is an analogous
question about $A_6$ requiring a more delicate argument since it
is nearly possible to acheive an extent.

In the language of section \ref{proofs}, the longest chain
generated by $P$ and $B$ has length 30.

This proof gives the idea of how to construct an extent: use odd
permutations for $Z$ and $W$ but even permutations for $P$ and
$B$.  Any method with these properties has a chance of working.
This idea leads to results of Saddleton (see theorems 4.8 and 4.11
of \cite{WhiteRing}).

In practice another type of lead, called a single lead,
is used to obtain an extent.

\bigskip

\subsection{Grandsire}

\bigskip

We now consider another method, the last of this article.
The method named Grandsire
(pronounced grand-sir)
is rung on an odd number of bells.
It was developed in the 1650's by Robert Roan on 5 bells,
and extensions to 7 and more bells took place in the
late 1600's or later.
The problem we referred to in
the abstract is on 7 bells, but first we explain
Grandsire on 5 bells.

The hunting subgroup $H_{10}$ is generated by $X=(1\ 2)(3\ 4)$ and
$Y=(2\ 3)(4\ 5)$ as usual. We introduce $Z=(1\ 2)(4\ 5)$, but the
first difference in Grandsire from Plain Bob is that we do $Z$ at
the very start.  This is irrelevant from a mathematical point of
view. Then we do $Y$, and then alternate $X$ and $Y$ until we have
run through the coset $ZH_{10}$. The last permutation done will be
$Y$, and in total we will have done $ZYXYXYXYXY=ZX$. Then we
repeat the permutations, i.e., do $Z, Y, X, Y, \ldots , X, Y$
until we have run through the coset $(ZXZ)H_{10}$, and then we
repeat the permutations again, running through $(ZXZXZ)H_{10}$,
and then we are back at rounds. Here are the rows of 3 plain
leads, a plain course. Neglecting the first row which is the
initial rounds, each column is a plain lead.

\bigskip

\begin{tabular}{rlllrl}
 1\ 2\ 3\ 4\ 5 & &\\
 & $Z$ &          & $Z$ & & $Z$ \\
 2\ 1\ 3\ 5\ 4 & & 2\ 1\ 5\ 4\ 3 & & 2\ 1\ 4\ 3\ 5\\
 & $Y$ &          & $Y$ & & $Y$ \\
 2\ 3\ 1\ 4\ 5 & & 2\ 5\ 1\ 3\ 4 & & 2\ 4\ 1\ 5\ 3 \\
 & $X$ &          & $X$ & & $X$ \\
 3\ 2\ 4\ 1\ 5 & & 5\ 2\ 3\ 1\ 4 & & 4\ 2\ 5\ 1\ 3\\
 & $Y$ &          & $Y$ & & $Y$ \\
 3\ 4\ 2\ 5\ 1 & & 5\ 3\ 2\ 4\ 1 & & 4\ 5\ 2\ 3\ 1\\
 & $X$ &          & $X$ & & $X$ \\
 4\ 3\ 5\ 2\ 1 & & 3\ 5\ 4\ 2\ 1 & & 5\ 4\ 3\ 2\ 1\\
 & $Y$ &          & $Y$ & & $Y$ \\
 4\ 5\ 3\ 1\ 2 & & 3\ 4\ 5\ 1\ 2 & & 5\ 3\ 4\ 1\ 2\\
 & $X$ &          & $X$ & & $X$ \\
 5\ 4\ 1\ 3\ 2 & & 4\ 3\ 1\ 5\ 2 & & 3\ 5\ 1\ 4\ 2\\
 & $Y$ &          & $Y$ & & $Y$ \\
 5\ 1\ 4\ 2\ 3 & & 4\ 1\ 3\ 2\ 5 & & 3\ 1\ 5\ 2\ 4\\
 & $X$ &          & $X$ & & $X$ \\
 1\ 5\ 2\ 4\ 3 & & 1\ 4\ 2\ 3\ 5 & & 1\ 3\ 2\ 5\ 4\\
 & $Y$ &          & $Y$ & & $Y$ \\
 1\ 2\ 5\ 3\ 4 & & 1\ 2\ 4\ 5\ 3 & & 1\ 2\ 3\ 4\ 5
\end{tabular}

\bigskip

It is because the first lead head is the result of $ZX=(1\ 2)(4\
5)(1\ 2)(3\ 4)=(3\ 4\ 5)$ which has order 3, that we get back to
rounds after $3\times 10=30$ permutations, and a plain course has
3 plain leads.

\bigskip

As with the Plain Bob method, we
sometimes add bob leads to the above plain course
to obtain a longer method.
The bob lead uses the permutation $Z$
applied instead of the last $X$ in a plain lead,
which would be
{\it two places before} the next appearance
of $Z$ in a plain course.

\begin{tabular}{rll}
 1\ 2\ 3\ 4\ 5 \\  & $Z$ &        \\  2\ 1\ 3\ 5\ 4 \\
& $Y$ &         \\  2\ 3\ 1\ 4\ 5 \\  & $X$ &         \\  3\ 2\ 4\
1\ 5 \\  & $Y$ &       \\  3\ 4\ 2\ 5\ 1 \\  & $X$ & \\ 4\ 3\ 5\
2\ 1 \\  & $Y$ &       \\  4\ 5\ 3\ 1\ 2 \\  & $X$ & \\
 5\ 4\ 1\ 3\ 2 \\  & $Y$ & \\  5\ 1\ 4\ 2\ 3 \\  & $Z$
&        \\  1\ 5\ 4\ 3\ 2 \\  & $Y$ & \\  1\ 4\ 5\ 2\ 3 \\ & $Z$&
\end{tabular}

\bigskip

If we use the bob before the first lead head, as shown above, the
first lead head will be $1\ 4\ 5\ 2\ 3$ instead of $1\ 2\ 5\ 3\
4$. This is the result of $(2\ 4)(3\ 5)$ applied to rounds. Since
this permutation has order 2, we will return to rounds after using
the bob twice in that place. In other words, the method consisting
of leads $B, P, P, B, P, P$ increases the number of permutations
from $30$ to $60$.

It is reasonable to ask, as usual,
if we could obtain a larger set
of permutations by using plain and bob leads in a
different arrangement.
The answer is no.  To see this, simply note that each
of $X, Y, Z$ is even.
Therefore the largest possible number of permutations
they can generate is 60 (the order of $A_5$), which
in fact is the case as we have shown.

\begin{theorem} \label{grandoubles}
There does not exist an extent of Grandsire on 5 bells
using plain and bob leads.
The longest possible method using plain and bob leads has $60$
permutations, and there does exist such a method.
\end{theorem}

To obtain the maximum number
of permutations on 5 bells
we would need to use an odd permutation,
which involves another type of lead called a
single lead.  We do not discuss single leads in
this article.

\bigskip

Next we consider Grandsire on 7 bells, or Grandsire Triples as it
is known to ringers, which is more interesting mathematically than
Grandsire on 5 bells. In this case $X=(1\ 2)(3\ 4)(5\ 6)$ and
$Y=(2\ 3)(4\ 5)(6\ 7)$ generate the hunting subgroup $H_{14}$, and
Grandsire Triples uses $Z=(1\ 2)(4\ 5)(6\ 7)$. As on 5 bells, we
do $Z$ first and then alternate $Y$ and $X$, and repeat. The first
lead head is the result of $ZX=(1\ 2)(4\ 5)(6\ 7)(1\ 2)(3\ 4)(5\
6)=(3\ 4\ 6\ 7\ 5)$ on rounds, which is $1\ 2\ 5\ 3\ 7\ 4\ 6$.
Since $ZX$ has order 5 there are 5 plain leads and 70 permutations
in a plain course of Grandsire on 7 bells.

To extend the method we  use  bob leads. The bob lead uses the
permutation $Z$ applied instead of the last $X$ in a plain lead,
as in Grandsire on 5 bells. If we do a bob lead on the earliest
possible occasion, the first lead head would become $1\ 7\ 5\ 2\
6\ 3\ 4$. This is the result of $B=(2\ 4\ 7)(3\ 6\ 5)$ acting on
rounds.  Since $B$ has order 3, we obtain a total of 210
permutations if we use $B, P, P, P, P, B, P, P, P, P, B, P, P, P,
P$.

It is possible to obtain  larger sets of permutations by using
plain and bob leads in different sequences. We may then reasonably
wonder as to the largest set we can get.  On 5 bells we used that
fact that $X, Y, Z$ are even to obtain an upper bound (which was
60). This argument will not work here, since $X$ and $Y$ are odd.
It is, in fact, conceivable that we could achieve an extent of all
$7!=5040$ permutations. As before, it is enough to consider the
action of $P$ and $B$ on lead heads. In 1741 John Holt came very
close to an extent and obtained $4998$ permutations.  This gave
rise to the famous old question mentioned in the introduction and
section \ref{groupform}:

\bigskip

{\bf Famous Old Question}.  Is it possible to ring all the $5040$
permutations on seven bells using the Grandsire method and plain
and bob leads only? In other words, does there exist an extent of
Grandsire on 7 bells using plain and bob leads?

\bigskip

Considering only lead heads, first note that there would be
$5040/14=360$ lead heads. Next check that $P=(3\ 4\ 6\ 7\ 5)$ and
$B=(2\ 4\ 7)(3\ 6\ 5)$ generate $A_6$. If $P$ and $B$ generate
$A_6$ unicursally then we would obtain the extent we are looking
for. We therefore arrive at the following question in order to
answer the Famous Old Question:

\bigskip

{\bf Question A}. Is $A_6$ (acting on $\{2,3,4,5,6,7\}$) generated
unicursally by $P=(3\ 4\ 6\ 7\ 5)$ and $B=(2\ 4\ 7)(3\ 6\ 5)$ ?

\subsection{A bob for the new method}

We must define a bob lead for the "new" method proposed in section
\ref{ccc}.  We propose the permutation $W=(3\ 4)$ instead of the
last $Z$ in a lead.  Thus a plain lead has the permutations
\[
X,Y,X,Y,X,Y,X,Y,X,Z
\]
and a bob lead has the permutations
\[
X,Y,X,Y,X,Y,X,Y,X,W.
\]
The sequence of leads $P,P,P,B,P,P,P,B,P,P,P,B$ yields an extent
of all 120 permutations.

\bigskip
\section{Proofs} \label{proofs}
\bigskip

We now proceed to give Rankin's proof that question A has a
negative answer. As we explained in section \ref{leads}, this
implies that the answer to the Famous Old Question is no. The
proof is by contradiction.

\medskip
We begin by supposing that $A_6$ is
generated unicursally by
$P=(3\ 4\ 6\ 7\ 5)$ and $B=(2\ 4\ 7)(3\ 6\ 5)$.
Assume there exists an ordering $g_1, g_2, \ldots ,g_{360}$
of the elements of $A_6$,
with the property that for each $i$,
$g_{i+1}=g_{i}t_i$ for some $t_i\in T=\{P, B\}$
(with subscripts modulo 360).
We will call any sequence of elements
$g_1, \ldots , g_m$ a {\it chain} of length $m$
if it has the property that for each $i$,
$g_{i+1}=g_{i}t_i$ for some $t_i\in T=\{P, B\}$.
With subscripts modulo $m$,
we consider a chain to be an infinite cyclic sequence.
If $g_{i+1}=g_{i}P$ then we say that $g_i$ is
{\it acted on} by $P$, and similarly for $B$.
Each element of $A_6$ is acted on by exactly
one of $P$ and $B$, by assumption.

The key is to consider left cosets of the cyclic
subgroup $C$ of order 5 generated by
$$\gamma= BP^{-1}=(2\ 4\ 7)(3\ 6\ 5)(3\ 5\ 7\ 6\ 4)=
(2\ 3\ 4\ 6\ 7).$$
This was noted by Thompson, who called these
cosets ``Q-sets".
We present Rankin's argument in a series of
observations, each one a lemma.
The idea of the proof is to show that, under a certain
transformation, the parity of the
number of chains remains constant.

\begin{lemma}\label{lemma1}
Every element of a coset $xC$ of $C$ is acted on by the same
element.
\end{lemma}

Proof: Suppose $x\gamma^i$ is acted on by $P$ (a similar argument
holds for $B$). Then the next element in the chain is
$x\gamma^iP$. But
$x\gamma^iP=x\gamma^{i-1}(BP^{-1})P=x\gamma^{i-1}B$. To avoid
repetition therefore, $B$ cannot act on $x\gamma^{i-1}$, so $P$
must act on $x\gamma^{i-1}$. This argument is valid for any $i$.

\bigskip

We now consider where the elements of a coset $xC$ appear in the
chain. For each $i$ between 1 and 5, we define a positive integer
$k_i$ between 1 and 5, by letting $x\gamma^{k_i}$ be the next
element of $xC$ in the chain after $x\gamma^{i}$. This defines a
permutation in $S_5$, which in two-line notation is
$$\sigma (x)=\pmatrix{1&2&3&4&5\cr
k_1&k_2&k_3&k_4&k_5\cr}$$ for each coset $xC$.

\begin{lemma}\label{lemma2}
The permutation $\sigma (x)$ is a 5-cycle.
\end{lemma}

Proof: Our assumption on the existence of a chain of length 360
implies that $\sigma (x)$ is a 5-cycle.

\bigskip

We will now rearrange the chain.  It is possible that
the `length 360 chain'
property may be destroyed during the rearrangement,
i.e., the single chain of length 360
may become several (disjoint) chains of smaller length.

Here is how the rearrangement is done. Let $xC$ be a coset acted
on by $B$ (if there is no such coset then every element of $A_6$
is acted on by $P$, which is impossible if there is only one
chain). Then the next element in the chain after $x\gamma^i$ is
$x\gamma^iB$. The chain can be divided up into 5 segments with
respect to $xC$, each segment being a sequence beginning with
$x\gamma^iB$ and ending with $x\gamma^{k_i}$. By definition of
$k_i$, there are no elements of $xC$ in the chain from
$x\gamma^iB$ to the element immediately preceding $x\gamma^{k_i}$.
In other words, in these segments with respect to $xC$, the only
element in a segment that is in $xC$ is the last element.
$$ \cdots x\gamma^i][x\gamma^iB \cdots x\gamma^{k_i}][x\gamma^{k_i}B \cdots $$
We permute these segments, so that the segment after $x\gamma^i$
now begins with $x\gamma^{i-1}B$ (and ends with
$x\gamma^{k_{i-1}}$).
$$ \cdots x\gamma^i][x\gamma^{i-1}B \cdots x\gamma^{k_{i-1}}]
[x\gamma^{k_{i-1}-1}B \cdots $$

\begin{lemma}\label{lemma3}
After the rearrangement, the coset $xC$ is acted on by $P$.  All
other cosets are unaffected, in terms of whether they are acted on
by $P$ or $B$.
\end{lemma}

Proof: Note that
$x\gamma^{i-1}B=x\gamma^{i-1}BP^{-1}P=x\gamma^iP$, so the element
after $x\gamma^i$ is $x\gamma^iP$. This implies that $x\gamma^i$,
and therefore the coset $xC$, is now acted on by $P$ in the
rearrangement. The proof of the second part is clear from the
construction.

\begin{lemma}\label{lemma4}
After the rearrangement there may be more than one chain, but
every chain contains an element of $xC$.
\end{lemma}

Proof: The rearrangement may alter the number of chains because,
for example, we may have $k_1=3$ and $k_2=2$. In this case, after
rearranging, we would have
$$[x\gamma B\cdots x\gamma^3][x\gamma^2B\cdots x\gamma^2]$$
which is a chain, and the other segments would form at least
one other chain.

If there were a chain not containing any element of $xC$ after the
rearrangement, this chain would not be affected by the
rearranging, and would therefore have existed before the
rearrangement.  But there was only one chain before the
rearrangement!

\bigskip

The next element of $xC$ after $x\gamma^i$ in the new arrangement
is $x\gamma^{k_{i-1}}$, so we define a permutation $\tau(x)$ in a
similar manner to $\sigma(x)$:
$$\tau (x)=\pmatrix{1&2&3&4&5\cr
k_5&k_1&k_2&k_3&k_4\cr}.$$

In the following, ``cycles'' means disjoint cycles including
1-cycles, as usual.

\begin{lemma}\label{lemma5}
The number of cycles in $\tau(x)$ is equal to the number of chains
after the rearrangement.
\end{lemma}

Proof: This follows from Lemma \ref{lemma4}.

\begin{lemma}\label{lemma6}
$(1\ 2\ 3\ 4\ 5)\tau(x)=\sigma(x)$.
\end{lemma}

Proof: This is straightforward.

\begin{lemma}\label{lemma7}
The number of cycles in $\tau(x)$ is odd.
\end{lemma}

Proof: Suppose $\tau (x)$ has $k$ cycles. The number of cycles in
$(1\ 5)\tau(x)$ is $k+1$ if 1 and 5 are in the same cycle in
$\tau(x)$, and $k-1$ otherwise.  Hence the number of cycles in
$(1\ 4)(1\ 5)\tau(x)$ has the same parity as the number of cycles
in $\tau(x)$. Similarly the number of cycles in $(1\ 2)(1\ 3)(1\
4)(1\ 5)\tau(x)$ has the same parity as the number of cycles in
$\tau(x)$. Lemmas \ref{lemma2} and \ref{lemma6} now give the
result.

\bigskip

We point out that this proof works since 5 is odd, so any odd
number could be used. If an even number is used instead of 5 then
the proof shows that the parity changes. This will be used in the
proof of theorem \ref{newthms4}.

Let $r_x$ denote the number of
chains after the rearrangement with respect to $xC$.
The following fact is all-important.

\begin{lemma}\label{lemma8}
$r_x$ is odd.
\end{lemma}

Proof: Combine Lemmas \ref{lemma5} and \ref{lemma7}.

We now rearrange again with respect to another coset $yC$ that is
acted on by $B$. (If there is no such coset, skip to the Famous
Old Theorem below.) Let $r_y$ denote the number of chains after
the rearrangement with respect to $yC$. Again we define
$k_i=k_i(y)$ and $\sigma (y)$ in a similar manner. Lemma
\ref{lemma2} becomes

\begin{lemma}\label{lemma2'}
The number of chains (before rearrangement with respect to $yC$)
containing elements of $yC$ is equal to the number of cycles in
$\sigma (y)$.
\end{lemma}

Lemma \ref{lemma3} shows that after rearranging, $yC$ is acted on
by $P$. Lemma \ref{lemma4} becomes

\begin{lemma}\label{lemma4'}
The number of chains not containing elements of $yC$ remains
constant.
\end{lemma}

Proof: Shown in proof of Lemma \ref{lemma4}.

We define $\tau(y)$ similarly to $\tau(x)$. Lemma \ref{lemma5}
becomes

\begin{lemma}\label{lemma5'}
$r_x - (\hbox{number of cycles in }\sigma (y))=
r_y - (\hbox{number of cycles in }\tau(y))$.
\end{lemma}

Proof: The lefthand side is the number of chains before
rearrangement not containing elements of $yC$, and the righthand
side is the number of such chains after rearrangement.

Lemma \ref{lemma6} remains the same, with $y$ in place of $x$, and
Lemma \ref{lemma7} becomes

\begin{lemma} \label{lemma7'}
$(\hbox{number of cycles in }\tau(y) ) \equiv
(\hbox{number of cycles in }\tau(x)) \ \hbox{ mod } \ 2$.
\end{lemma}

Proof: Same as Lemma \ref{lemma7}.

Again, the following will be important.

\begin{lemma} \label{lemma8'}
$r_y$ is odd.
\end{lemma}

Proof: By Lemma \ref{lemma8} and Lemmas \ref{lemma5'} and
\ref{lemma7'}.

\bigskip

The series of lemmas shows that after this second
rearrangement with respect to $yC$, the number of chains
is still odd.
Repeat the rearrangement with respect to
every coset acted on by $B$
until all cosets are acted on by $P$.
The point is that after each rearrangement with respect to a coset
acted on by $B$, the coset is now acted on by $P$, and
also we can apply the lemmas and conclude that the
number of chains remains odd.

\begin{theorem}  $($Thompson$)$
\begin{enumerate} \item  $A_6$
(acting on $\{2,3,4,5,6,7\}$) is not generated unicursally by
$P=(3\ 4\ 6\ 7\ 5)$ and $B=(2\ 4\ 7)(3\ 6\ 5)$. \item It is not
possible to ring all the $5040$ permutations on seven bells using
the Grandsire method and plain and bob leads. In other words,
there does not exist an extent of Grandsire on 7 bells using plain
and bob leads.  \end{enumerate}
\end{theorem}

Proof:  (Rankin) We showed in section \ref{leads} that (2) is
implied by (1).  To prove (1), by the above lemmas we may assume
that all cosets are acted on by $P$. Then every element of $A_6$
is acted on by $P$. The chains we have must therefore be the
cosets of the subgroup $M$ generated by $P$. By Lemma
\ref{lemma8'}, the number of chains is odd. However, since $P$ has
order 5, there are $|A_6:M|=360/5=72$ cosets, which is an even
number. This contradiction proves the theorem.

\bigskip

We remark that the roles of $P$ and $B$ could be interchanged in
the above proof, since the subgroup generated by $B$ also has even
index. 

The following is shown in \cite{Dick}.

\begin{corollary} The largest number of permutations that can be
rung on seven bells, using the Grandsire method and plain and bob
leads, is $4998$.
\end{corollary}

Proof: We know that one chain of length 360 (in $A_6$) is not
possible. The shortest possible chain has length 3 since $B$ has
order 3, so the longest possible chain has length $\leq 357$. A
chain with 357 leads has $357\times 14 =4998$ permutations.

As we said earlier, there does exist a method with 4998
permutations due to Holt, so 4998 is best possible.

\bigskip

\section{A New Result}\label{newthm}

\bigskip

We next prove a small result concerning problem 1 of section
\ref{groupform}. Observe that the result is not trivial since
there are $2^{24}$ possible orderings of $A$ and $B$.

\begin{theorem} \label{newthms4}
$S_4$ is not unicursally generated by $A=(3\ 4)$
and $B=(1\ 2\ 3)$.
\end{theorem}

Proof: Suppose to the contrary that $S_4$ is unicursally generated
by $A$ and $B$. As in the sequence of lemmas, we rearrange with
respect to left cosets of $C=\langle A^{-1}B\rangle$ that are
acted on by $A$. After the rearrangement, such a coset is acted on
by $B$.

Since $C$ has order 4 which is even, the proof of Lemma
\ref{lemma7} shows that after the rearrangement with respect to
$xC$, the number of cycles in $\tau(x)$ is even, and in general
that the parity of the number of cycles changes after each
rearrangement with respect to a coset of $C$.

Suppose we perform a total of $m$ rearrangements to get all cosets
acted on by $B$. By Lemma \ref{lemma3} second sentence, the number
of cosets acted on initially by $A$ is $m$.  Note that $m\leq 6$
as $|C|=4$. After all $m$ rearrangements, the number of chains has
the opposite parity to $m$, by the previous paragraph.

By assumption, we start with one chain of length 24.
If $m$ is even then
the number of chains (after all rearrangements)
is odd.  The proof of Rankin's theorem above shows that
$\langle B \rangle$ has odd index, a contradiction.

The only alternative is that $m$ is odd, and we have 1, 3, or 5
cosets of $C$ acted on by $A$ initially. If 1 coset was acted on
by $A$, then 5 cosets were acted on by $B$. Therefore 20 elements
of $S_4$ were acted on by $B$, which implies there are three
consecutive $B$'s. This is not possible as $B$ has order 3. If 5
cosets were acted on by $A$, then 1 coset was acted on by $B$, and
a similar argument gives a contradiction.

The final possibility is that 3 cosets were acted
on by $A$ and 3 by $B$.
Then 12 elements of $S_4$ are acted on by each
of $A$ and $B$.
Since $A$ has order 2, we must have $A$ and $B$
alternating in the chain.
But $AB$ and $BA$ have order 4, so neither of these
work.

\end{document}